\documentclass[12pt,english,fleqn,liststotoc,bibtotoc,idxtotoc,BCOR7.5mm,tablecaptionabove]{amsart}
\usepackage[T1]{fontenc}
\usepackage[latin1]{inputenc}
\usepackage{geometry}
\usepackage{color}
\usepackage{amsthm}
\usepackage{amstext}
\usepackage{amssymb}
\usepackage{amsmath}
\usepackage{graphicx}
\usepackage{young}
\usepackage[section]{placeins}
\usepackage[svgnames]{xcolor}
\usepackage{breakurl}
\usepackage{ifpdf}
\ifpdf

\IfFileExists{lmodern.sty}
{\usepackage{lmodern}}{}
\makeatletter

\usepackage[unicode=true,
bookmarks=true,bookmarksnumbered=true,bookmarksopen=false,
breaklinks=false,pdfborder={0 0 1},backref=false,colorlinks=true]
{hyperref}
\hypersetup{
	pdfauthor={Nick Early},
	linkcolor=black, citecolor=black, urlcolor=blue, filecolor=blue,  pdfpagelayout=OneColumn, pdfnewwindow=true,  pdfstartview=XYZ, plainpages=false, pdfpagelabels}

\theoremstyle{plain}
\newtheorem{thm}{\protect\theoremname}
\theoremstyle{plain}
\newtheorem{conjecture}[thm]{\protect\conjecturename}
\theoremstyle{plain}
\newtheorem{question}[thm]{\protect\questionname}
\theoremstyle{remark}

\theoremstyle{plain}

\theoremstyle{plain}

\theoremstyle{remark}

\theoremstyle{remark}

\theoremstyle{definition}

\theoremstyle{definition}
\newtheorem{example}[thm]{\protect\examplename}
\theoremstyle{plain}

\makeatother

\providecommand{\claimname}{\inputencoding{latin9}Claim}
\providecommand{\conjecturename}{\inputencoding{latin9}Conjecture}
\providecommand{\corollaryname}{\inputencoding{latin9}Corollary}
\providecommand{\definitionname}{\inputencoding{latin9}Definition}
\providecommand{\examplename}{\inputencoding{latin9}Example}
\providecommand{\lemmaname}{\inputencoding{latin9}Lemma}
\providecommand{\notename}{\inputencoding{latin9}Note}
\providecommand{\propositionname}{\inputencoding{latin9}Proposition}
\providecommand{\questionname}{\inputencoding{latin9}Question}
\providecommand{\remarkname}{\inputencoding{latin9}Remark}
\providecommand{\theoremname}{\inputencoding{latin9}Theorem}

\newcommand{\symm}{\mathfrak{S}}

\title[Conjectures for $h^*$-vectors of Dilated Hypersimplices]{Conjectures for Ehrhart $h^*$-vectors of Hypersimplices and Dilated Simplices}
\author{Nick Early}
\thanks{The author was partially supported by RTG grant NSF/DMS-1148634,\\
	University of Minnesota, email: \href{mailto:earlnick@gmail.com}{earlnick@gmail.com}}
\begin{document}

\maketitle

\begin{abstract}

We formulate conjectures giving combinatorial interpretations of the Ehrhart $h^*$-vector, for hypersimplices, for dilated simplices and for generic cross-sections of cubes, in terms of certain decorated ordered set partitions.   All were formulated and checked computationally during our graduate study at Penn State.

\end{abstract}
	\begingroup
\let\cleardoublepage\relax
\let\clearpage\relax
\tableofcontents
\endgroup

\section{Discussion}

Fix integers $1\le a<n$ and put $b=n-a$.  Let $(L_1,L_2,\ldots, L_k)$ be an ordered set partition of $\{1,\ldots, n\}$ and $(l_1,\ldots, l_k)$ be an ordered partition of $a$.  Denote by $B_{a,b}=\{x\in\lbrack0,1\rbrack^{n}:\sum x_i=a\}$, the hypersimplex given as the convex hull of all permutations of the vector $(1,\ldots, 1,0,\ldots, 0)$ with $a$ ones and $b$ zeros.  

Call the (decorated) ordered set partition $((L_1)_{l_1},\ldots, (L_k)_{l_k})$ \textit{hypersimplicial} if we have $1\le l_i <\vert L_i\vert$ for each $i=1,\ldots, k$.  From Theorem 5.1 of Ocneanu's preprint \cite{OcneanuPermutations}, it follows that the hypersimplicial ordered set partitions are enumerated (up to cyclic rotation of blocks) by the Eulerian numbers.  

While the conjecture here is stated in almost entirely combinatorial terms, we remark that stronger, geometric and representation theoretic results are known, the most immediate of which is that the dimension of the $\mathfrak{S}_n$-module associated to plates in a hypersimplex $B_{a,b}$ is the same as the relative volume of that hypersimplex, the Eulerian number $E_{a-1,b-1}$, which counts the number of permutations in $\mathfrak{S}_{n-1}$ having $a-1$ descents.  This dimension count can be seen from geometric arguments for counting plates in integer dilations of a simplex, together with an application of the classical Worpitzky identity, as mentioned below Theorem 9 in \cite{EarlyPlateAnnouncement}, to which we refer for further discussion of terminology and context.
\begin{figure}[h!]
	\centering
	\includegraphics[width=.55\linewidth]{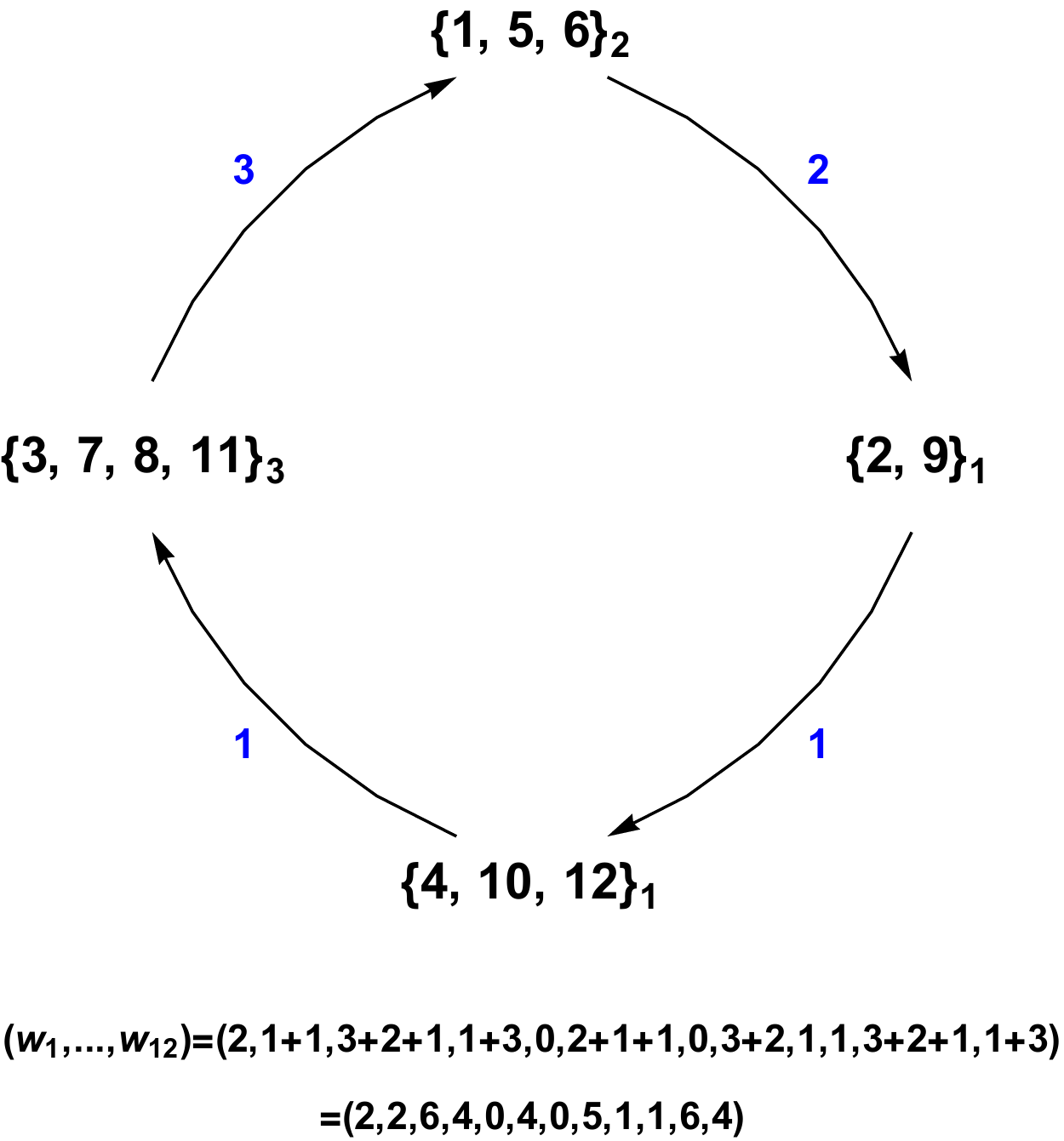}
	\caption{The winding number of $ (\{1,5,6\}_2,\{2,9\}_1,\{4, 10, 12\}_1,\{3,7,8,11\}_3)$ is $35/7=5$.}
	\label{fig:clock-example}
\end{figure}

We introduce two labeling schemes which extend to decorated ordered set partitions more generally.  First, for each decorated ordered set partition $((L_1)_{l_1},\ldots, (L_k)_{l_k})$, we define an $n$-tuple $(p_1,\ldots, p_n)\in (\mathbb{Z}\slash a)^n$, where $p_j=l_1+\cdots+l_{j-1}$ is the cumulative position of the block containing $j$ from the first block.  Thus, in particular, for any $i\in L_1$, we have $p_i=0$.  

Second, fix the cyclic order $\sigma=(12\cdots n)$ and a decorated ordered set partition 
$$P=((L_1)_{l_1},\ldots, (L_k)_{l_k})$$
with $1\in L_1$.  Denote by $w_{i}$ the nonnegative representative of $p_{i+1}-p_{i}$ modulo $a$ in $\{0,\ldots, a-1\}$, where the indices $i$ are taken modulo $n$.  One can check that the \textit{winding map} $P\mapsto (w_1,\ldots, w_n)$ is one-to-one. Call $\sum_{i=1}^n w_i$ the \textit{level} and $w_P=\frac{1}{a}\sum_{i=1}^n w_i$ the winding number of $P$.  Note that $(w_1,\ldots, w_n)$ could also be obtained directly from $P$, as illustrated in Figure \ref{fig:clock-example}.

\section{$h^*$-vectors of hypersimplices}

\begin{conjecture}
	Denote by $m_0,m_1,\ldots$ the numbers of hypersimplicial ordered set partitions $((L_1)_{l_1},\ldots, (L_k)_{l_k})$ for the hypersimplex $B_{a,b}$, with $1\in L_1$, having winding numbers respectively $w_P=0,1,2,\ldots$.  Then, the vector $(m_0,m_1,\ldots)$ equals the Ehrhart $h^*$-vector for the hypersimplex $B_{a,b}$.
\end{conjecture}

See for example \cite{Braun} for discussion about the $h^*$-vector.

\begin{figure}[h!]
	\centering
	\includegraphics[width=1\linewidth]{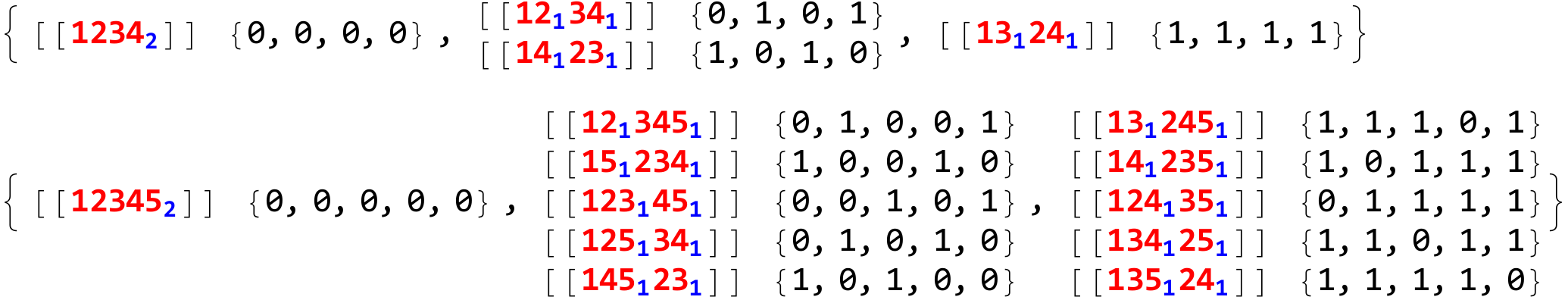}
	\caption{The column counts are the $h^*$-vectors respectively $(1,2,1)$ and $(1,5,5)$ for $B_{2,2}$ and $B_{2,3}$.}
	\label{fig:hypersimplex-conjecture-examples}
\end{figure}

\begin{example}
	We illustrate in Figure \ref{fig:hypersimplex-conjecture-examples} the labeling and enumeration related to the hypersimplices $B_{2,2}$ and $B_{2,3}$.
\end{example}

\begin{question}
	Is there a correspondence between hypersimplicial ordered set partitions and the unimodular simplices in the known triangulation of the hypersimplex?
\end{question}

\section{$h^*$-vectors of dilated simplices}
One may relax the hypersimplex condition on the possible values of the decorations of blocks for the decorated set partitions.  Say that a (decorated) ordered set partition $((L_1)_{l_1},\ldots, (L_k)_{l_k})$ is \textit{simplicial} if the $l_i$ satisfy $1\le l_i\le r-1$ and $\sum l_i=r$, for the $(n-1)$-dimensional simplex 
$$\Delta_r^n=\left\{x\in\lbrack 0,r\rbrack^n: \sum x_i=r \right\}.$$
which has been dilated by an integer factor $r$.  One finds, as announced in \cite{EarlyPlateAnnouncement}, that there are $r^{n-1}$ simplicial ordered set partitions $((L_1)_{l_1},\ldots, (L_k)_{l_k})$ with $1\in L_1$.

\begin{figure}[h!]
	\centering
	\includegraphics[width=1\linewidth]{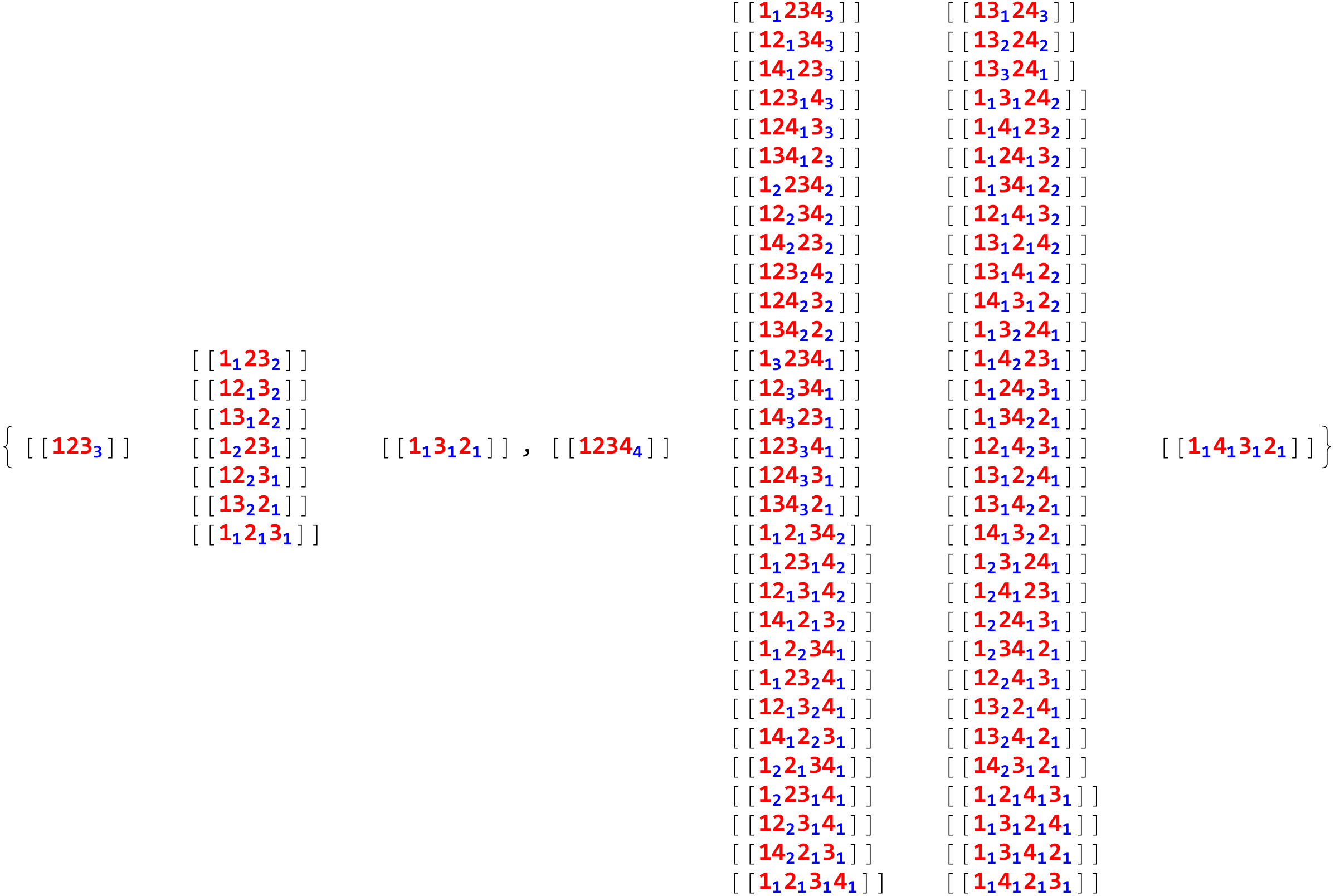}
	\caption{Left: decomposition of simplicial ordered set partitions with $1\in L_1$, for $\Delta_3^3$.  Right: same for $\Delta_4^4$.}
	\label{fig:simplex-plate-module-dilation-graded}
\end{figure}

In Figure \ref{fig:simplex-plate-module-dilation-graded}, simplicial ordered set partitions are listed for $\Delta_3^3$ and $\Delta_4^4$, decomposed according to winding number. The winding map from the Discussion is now a bijection between simplicial ordered set partitions for $\Delta_r^n$ and the set
$$\left\{x\in\left(\mathbb{Z}\slash r\right)^n:\sum x_i\equiv 0\text{ mod}(r)\right\}.$$

Generating functions for the numbers of the simplicial ordered set partitions are given below.  The coefficients of the numerator in the table in position $(r,n)$ count the number of simplicial ordered set partitions for $\Delta_r^n$, when graded by winding number.  For instance, for $\Delta_3^3$ there are respectively $1,7$ and $1$ simplicial ordered set partitions with winding numbers $0,1,2$.
$$
\begin{array}{cccc}
r\backslash n & 2 & 3 & 4 \\
2&\frac{1+x}{(1-x)^2} & \frac{1+3 x}{(1-x)^3} & \frac{1+6x+x^2}{(1-x)^4} \\
3&\frac{1+2 x}{(1-x)^2} & \frac{1+7 x+x^2}{(1-x)^3} & \frac{1+16x+10 x^2}{(1-x)^4} \\
4&\frac{1+3 x}{(1-x)^2} & \frac{1+12x+3 x^2}{(1-x)^3} & \frac{1+31 x+31 x^2+x^3}{(1-x)^4} \\
\end{array}
$$

It is easy to see that that the Ehrhart polynomial of $\Delta_r^n$, which counts the number of integer lattice points in its $s^\text{th}$ dilate, is $\binom{n-1+r s}{n-1}$.  Then the $h^*$-vector is the numerator of the generating function for the series
\begin{eqnarray}\label{eqn: dilatedSimplexSeriesCoefficients}
\sum_{s=0}^\infty \binom{n-1+r s}{n-1}x^s.
\end{eqnarray}

\begin{conjecture}
		Denote by $m_0,m_1,\ldots$ the numbers of simplicial ordered set partitions $$((L_1)_{l_1},\ldots, (L_k)_{l_k})$$ for the dilated simplex $\Delta_r^n$, with $1\in L_1$, having winding numbers respectively $w_P=0,1,2,\ldots$.  Then, the vector $(m_0,m_1,\ldots)$ equals the Ehrhart $h^*$-vector for the simplex $\Delta_r^n$, and can be obtained by summing the series in Equation \ref{eqn: dilatedSimplexSeriesCoefficients} for $r,n$ given.
\end{conjecture}
From the Conjecture it follows that the coefficients of the $h^*$-vector for $\Delta_r^n$ sum to $r^{n-1}$.
For example, for the tetrahedron $\Delta_4^4$, with $(r,n)=(4,4)$, we have
$$\sum_{s=0}^\infty\binom{4-1+4s}{4-1}x^s=\frac{1+31x+31x^2+x^3}{(1-x)^4},$$
where $1+31+31+1=4^{4-1},$ and where the coefficients correspond to the right four columns of Figure \ref{fig:simplex-plate-module-dilation-graded}.

\begin{figure}[h!]
	\centering
	\includegraphics[width=0.4\linewidth]{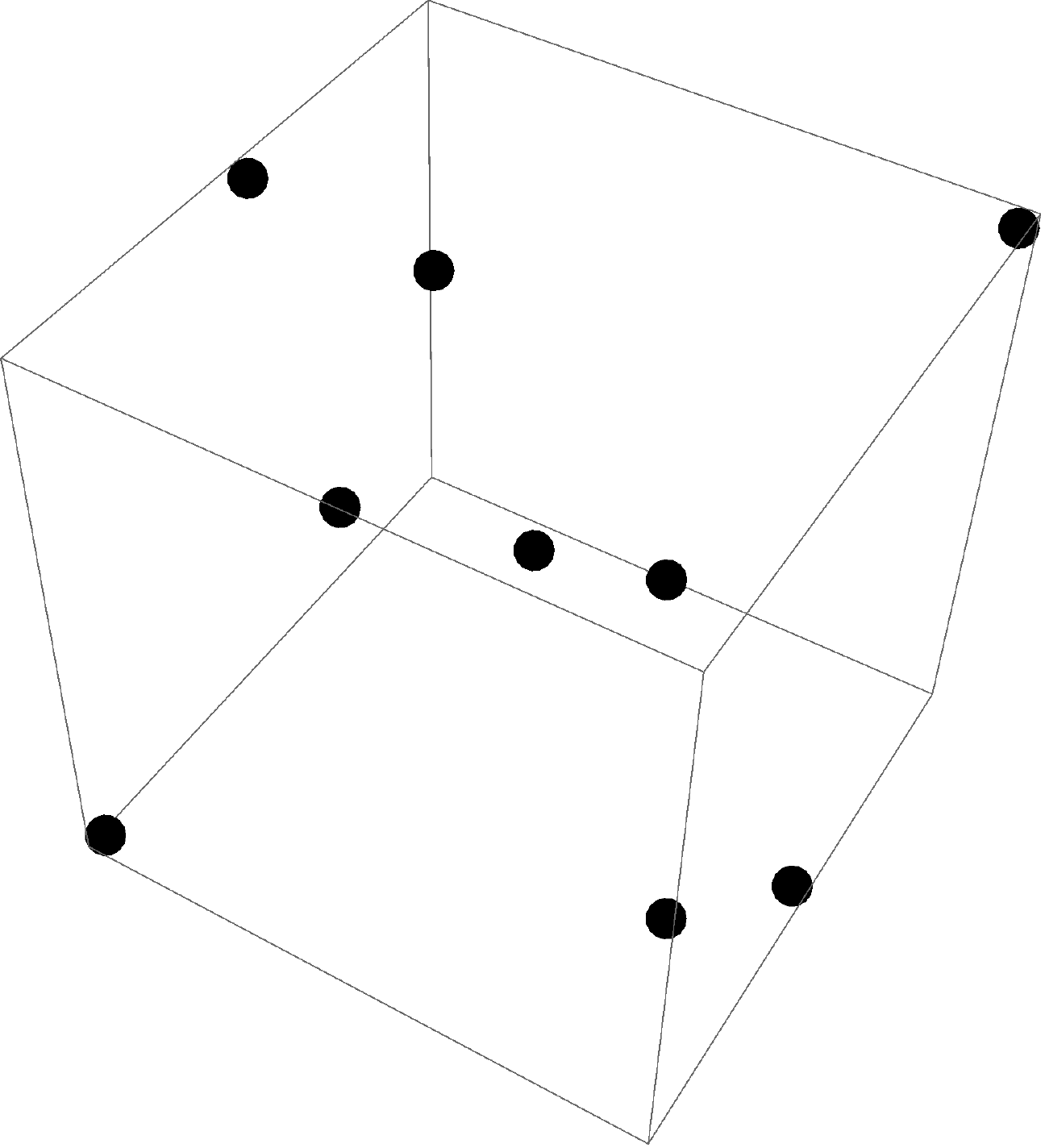}
	\caption{Decomposition into modular sections, by winding number, for simplicial ordered set partitions for $\Delta_3^3$.}
	\label{fig:winding-number-layers-3-variable-dilated}
\end{figure}
In Figure \ref{fig:winding-number-layers-3-variable-dilated}, the black dots label simplicial ordered set partitions in $\Delta_3^3$ in the left three columns of Figure \ref{fig:simplex-plate-module-dilation-graded}, and are located at the points
\begin{eqnarray*}
& (0,0,0) &\\
&(0,1,2),(0,2,1),(1,0,2),(1,1,1),(1,2,0),(2,0,1),(2,1,0)&\\
&(2,2,2)&
\end{eqnarray*}
all of which have residue zero modulo $3$ and winding numbers respectively $(0,3,6)/3=(0,1,2)$.

The modular cross-sections of $(\mathbb{Z}\slash r)^n$ are themselves integer lattice points of generalized hypersimplices!  In particular, as depicted in Figure \ref{fig:winding-number-layers-3-variable-dilated}, they are invariant under the full action of $\symm_n$, motivating in part an additional question.
\begin{question}
The decompositions of hypersimplicial and simplicial ordered set partitions are invariant only under $\mathbb{Z}\slash n$ and depend on the choice of an $n$-cycle.  Are there objects, corresponding to respectively hypersimplices and dilated simplices, which are preserved by $\symm_n$ and which are counted by the entries of the $h^*$-vector?
\end{question}

We have some encouraging initial observations in low dimensions.  We leave this to future work.

\section{$h^*$-vectors of generic cross-sections of cubes}

For each $s=1,\ldots, rn-1$, let $I_{r,s}^n=\left\{\lbrack0,r\rbrack^n: \sum_{i=1}^n x_i=s\right\}$.  Note that fixing $r$ and $n$ while letting $s$ vary gives all integer-sum hyperplane cross-sections of $\lbrack 0,r\rbrack^n$ perpendicular to the diagonal $(1,1,\ldots, 1)$.  

\begin{conjecture}
	Denote by $m_0,m_1,\ldots, $ the set of all decorated ordered set partitions $$((L_1)_{l_1},\ldots, (L_k)_{l_k})$$ satisfying $l_1+\cdots+l_k=s$ and $1\le l_i\le r\vert L_i\vert-1$, with $1\in L_1$, having winding numbers respectively $w_P=0,1,2,\ldots$.  Then, the vector $(m_0,m_1,\ldots)$ equals the $h^*$-vector for the hyperplane cross-section $I_{r,s}^n$.  
\end{conjecture} 

\begin{figure}[h!]
	\centering
	\includegraphics[width=1\linewidth]{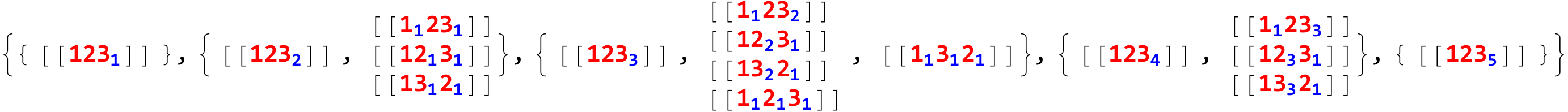}
	\caption{Combinatorial interpretation of $h^*$-vectors for cross-sections of $\lbrack 0,2\rbrack^3$}
	\label{fig:hyperplanecrosssectionsdim3}
\end{figure}

For example, for $n=3$ and $r=2$ we have the $h^*$-vectors
$$\{(1),(1,3),(1,4,1),(1,3),(1)\}$$
which decompose the relative volumes $(1,4,6,4,1)$ and which count the columns in Figure \ref{fig:hyperplanecrosssectionsdim3}.


\vspace{.5in}
\section{Acknowledgements}
We are grateful to Adrian Ocneanu for many intensive discussions on related topics throughout our graduate work.  We thank Alex Postnikov and Victor Reiner for useful comments and encouragement.


\end{document}